\documentclass[entropy,article,accept,moreauthors,pdftex,10pt,a4paper]{mdpi}

\usepackage{amssymb}
\usepackage{amsmath}
\usepackage{commath}
\usepackage[scr]{rsfso}
\usepackage{mathtools}
\usepackage{microtype}

\renewcommand{\[}{\begin{equation}} 
\renewcommand{\]}{\end{equation}} 
\def\<#1\>{\begin{align}#1\end{align}} 

\let\textcite\citet

\renewcommand{\v}{\relax\ifmmode\expandafter\boldsymbol\else\expandafter\textv\fi} 
\newcommand{\given}{\mid} 
\newcommand{\iid}[1][]{\overset{\textrm{iid}}{\sim}\ifx\\#1\\\else\operatorname{#1}\fi} 
\newcommand{\avtdist}[1][]{\sim\ifx\\#1\\\else\operatorname{#1}\fi} 
\newcommand{\R}{\mathbb{R}} 
\newcommand{\eps}{\varepsilon} 
\newcommand{\f}{\expandafter\operatorname} 
\newcommand{\I}[1][]{\operatorname{I}\ifx\\#1\\\else_{\{#1\}}\fi} 
\newcommand{\goesto}{\rightarrow} 

\newtheorem{conjecture}[theorem]{Conjecture}

\firstpage{1}
\articlenumber{391}
\doinum{10.3390/e19080391}
\pubvolume{19}
\pubyear{2017}
\copyrightyear{2017}
\history{Received: 20 June 2017; Accepted: 28 July 2017; Published: 29 July 2017}
\pdfoutput=1

\Title{A Noninformative Prior on a Space of Distribution~Functions}
\Author{Alexander Terenin and David Draper *}
\AuthorNames{Alexander Terenin and David Draper}
\address[1]{%
Applied Mathematics and Statistics, University of California, Santa Cruz,  CA 95064, USA; aterenin@ucsc.edu\\}
\corres{Correspondence: draper@ucsc.edu; Tel.: +1-831-459-1295
}

\abstract{In a given problem, the Bayesian statistical paradigm requires the specification of a prior distribution that quantifies relevant information about the unknowns of main interest external to the data.
In cases where little such information is available, the problem under study may possess an invariance under a transformation group that encodes a lack of information, leading to a unique prior---this idea was explored at length by E.T. Jaynes.
Previous successful examples have included location-scale invariance under linear transformation, multiplicative invariance of the rate at which events in a counting process are observed, and the derivation of the Haldane prior for a Bernoulli success probability.
In this paper we show that this method can be extended, by generalizing Jaynes, in two ways: (1) to yield families of approximately invariant priors; and (2) to the infinite-dimensional setting, yielding families of priors on spaces of distribution functions.
Our results can be used to describe conditions under which a particular Dirichlet Process posterior arises from an optimal Bayesian analysis, in the sense that invariances in the prior and likelihood lead to one and only one posterior distribution.
}

\keyword{Bayesian nonparametrics; Dirichlet process; functional equations; Hyers--Ulam--Rassias stability; improper prior; invariance; optimal Bayesian analysis; transformation group}

\begin{document}

\section{Introduction} \label{introduction}

Consider a statistician working on a problem $P$ in which a vector $\v{y} = ( y_1, .., y_n )$ of real-valued outcomes is to be observed, and---prior to, i.e., without observing $\v{y}$---the statistician's uncertainty is exchangeable, in the usual sense of being invariant under permutation of the order in which the outcomes are listed in $\v{y}$.  
This situation has extremely broad real-world applicability, including (but~not limited to) the analysis of a completely randomized controlled trial, in which participants---ideally, similar to elements of a population to which it is desired to generalize inferentially---are randomized.
Each participant is assigned either to a control group that receives the current best treatment, or to an experimental group that receives a new treatment whose causal effect on $y$ is of interest.
This design, while extremely simple, has proven to be highly useful over the past 90 years, in fields as disparate as agriculture \cite{fisher25}, medicine \cite{mrc48}, and (in contemporary usage) $A/B$ testing in data science on a massive scale~\cite{longbotham15}.
We use randomized controlled trials as a motivating example below, but we emphasize that they constitute only one of many settings to which the results of this paper apply.

Focusing just on the experimental group in the randomized controlled trial, the exchangeability inherent in $\v{y}$ implies via de Finetti's Theorem \cite{definetti37} that the statistician's state of information may be represented by the hierarchical model
\< \label{de-finetti-representation}
y_i \given F & \iid F
&
F & \avtdist\pi(F)
\>
for $i=1,..,n$, where $F$ is a cumulative distribution function (CDF) on $\R$ and $\pi(F)$ is a prior on the space of all such CDFs, i.e., the infinite-dimensional probability simplex $S_\infty$.
Note that (\ref{de-finetti-representation}) has uniquely specified the likelihood in a Bayesian nonparametric model for $\v{y}$, and all that remains is specification of $\pi(F)$.

Speaking \textls[-15]{now more generally (not just in the context of a randomized controlled trial), suppose that the nature of the problem $P$ enables the analyst to identify an alternative statistical problem $\tilde{P}$ in} which
\<
\tilde{P} &= g(P)
&
&\text{such that}
&
g &\in G \, ,
\>
where $G$ is a collection of transformations $g$ from one problem to another
having the property that, without having seen any data, $\tilde{P}$ and $P$ are the \emph{exact same problem}.
Then, the prior $\tilde{\pi}$ under $\tilde{P}$ must be the same as the prior $\pi$ under $P$!
Furthermore, since this holds for any $g \in G$, the result will be, as long as $G$ is endowed with enough structure, that there is one and only one prior $\pi$, for use in $P$, that respects the inherent invariance of the problem under study.
Bayes' Rule then implies that there is one and only one posterior distribution under $P$.
When this occurs---when both the likelihood function and the prior are uniquely specified, as in the example above---we say that the problem $P$ admits an \emph{optimal Bayesian analysis}.

The logic \textls[-15]{underlying the above argument has been used to motivate and formalize the notion of noninformative priors for decades.
Indeed, in the special case where $F$ is parametric and $G$ is a group of transformations encoding invariance with respect to monotonically-transformed units of measurement, \textcite{jeffreys46} derived the resulting prior distribution.
As another example, \textcite{jaynes03} derived the prior distribution for the mean number of arrivals of a Poisson process by using its characterization as a L\'evy counting process to specify an appropriate transformation group.
Notably, the resulting prior distribution} is \emph{not} the Jeffreys prior, because the problem's invariance and corresponding transformation group are different.
See \textcite{eaton89} for additional work on this subject.

Having studied this line of reasoning, it is natural to ponder its generality.
In this paper we show that the argument can be made quite general---we prove that the argument's formal notions
\begin{itemize}[leftmargin=2.3em,labelsep=4mm]
\item[(a)] can be generalized to include \emph{approximately} invariant priors in an $\eps$--$\delta$ sense; and
\item[(b)] can be extended to infinite-dimensional priors on spaces of CDFs.
\end{itemize}

We focus on the setting described in (\ref{de-finetti-representation}) and defer more general situations to future work.
In this setting we derive a number of results, ultimately showing that the Dirichlet Process \cite{ferguson73} prior $\f{DP}(\eps, F_0)$ is an approximately invariant stochastic process for any CDF $F_0$ on $\R$ and sufficiently small $\eps > 0$.
Together with de Finetti's Theorem, this demonstrates that the posterior distribution
\[
F \given \v{y} \avtdist[DP]\del{n, \hat{F}_n} \, ,
\]
where $\hat{F}_n$ is the empirical CDF, corresponds in a certain sense to an optimal Bayesian analysis---see Section \ref{discussion} for more on this point.

Not all approaches to noninformative priors are based on group invariance.
Perhaps the earliest approach can be traced back to \textcite{laplace74}, who proposed a Principle of Indifference under which, if all that is known about a quantity $\theta$ is that $\theta \in \Theta$ (for some set $\Theta$ of possible values), then the prior should be uniform on $\Theta$.
For example, consider $\Theta = (0,1)$: the fact that $\theta \avtdist[U](0,1)$ is not consistent with $f(\theta) \avtdist[U](0,1)$ for any monotonic nonlinear $f$ requires that the problem $P$ under study must uniquely identify the scale on which uniformity should hold for the principle to be valid---this was a major reason for the rise of non-Bayesian theories of inference in the 19th century \cite{hald07}.
\textcite{bernardo79} has proposed a notion of noninformative priors that is defined by studying their effect on posterior distributions, and choosing priors that ensure that prior impact is minimized.
\textcite{jaynes68} has proposed the Maximum Entropy Principle, which defines noninformative prior distributions via information-theoretic arguments, for use in settings in which invariance considerations do not lead to a unique prior.
All of these notions are different, and applicable to problems where the corresponding notions of noninformativeness arise most naturally.

Most of the work on noninformative priors has focused on the parametric setting, in which the number of unknown quantities is finite.~In contrast, Bush et al. \cite{bush10} and Lee et al. \cite{lee14} have derived results on noninformative priors in Dirichlet Process Mixture models.~Their notion of noninformativeness is completely different from our own, as it is a posteriori, i.e., it involves examining the behavior of the posterior distribution under the priors studied.
This makes their approach largely complementary to ours: in specifying priors, it is helpful to understand both the prior's effect on the posterior and the prior's behavior a priori without considering any data.

Here we study noninformative prior specification from a \emph{strictly a priori} perspective.
We do not consider the prior's effect on the posterior distribution.
There is no data or discussion of computation.

Our motivation is a generalization of the following argument by \textcite{jaynes68}.
Suppose that in the randomized controlled trial described above, the outcome $y$ of interest is binary. By de Finetti's Theorem, we know that
\[
y_i \given \theta_1 \iid[Ber](\theta_1)
\]
is the unique likelihood for (e.g., the treatment group in) this problem.
Suppose further that the statistician's state of information about $\theta_1$ external to the data set $\v{ y }$ is what Jaynes calls ``complete~initial ignorance'' except for the fact that $\v{ \theta } = ( \theta_1, \theta_2 )$ is such that
\[
\{ 0 \le \theta_1 \le 1, 0 \le \theta_2 \le 1, \theta_1 + \theta_2 = 1 \}
\,.
\]

Jaynes argues that this state of information is equivalent to the statistician possessing complete initial ignorance about all possible rescaled and renormalized versions of $\v{ \theta }$, namely
\[
\v{ \theta' } = \left( \frac{ c_1 \theta_1 }{ c_1 \theta_1 + c_2 \theta_2 }, \frac{ c_2 \theta_2 }{ c_1 \theta_1 + c_2 \theta_2 }\right)
\]
for all positive $c_1,c_2$.
Jaynes shows that this leads uniquely to the Haldane prior
\<
\pi( \theta_1 ) &\propto \frac{ 1 }{ \theta_1 ( 1 - \theta_1 ) }
&
&\text{or equivalently}
&
\pi(\theta_1, \theta_2) &\propto \frac{1}{\theta_1\theta_2} \, ,
\>
where $\theta_2 = 1 - \theta_1$.
Combining this result with the unique Bernoulli likelihood under exchangeability, in our language Jaynes has therefore identified an instance of optimal Bayesian analysis.
In what follows we (a) extend Jaynes's argument to the multinomial setting with $p$ outcome categories for arbitrary finite $p$ and (b) show how this generalization leads to a unique noninformative prior on $S_\infty$.

\begin{table}[H]
\begin{center}
\begin{tabular}{cl}
\toprule
\textbf{Expression} & \textbf{Description}\\\midrule
$\f{Dir}(\v\alpha)$ & The Dirichlet distribution with concentration vector $\v\alpha$.
\\
$\f{Dir}(\alpha, \v{F}_0)$ & The Dirichlet distribution with concentration parameter $\alpha$ and mean probability vector $\v{F}_0$.
\\
$\f{Dir}(0)$ & The improper Dirichlet distribution corresponding to $\lim_{\alpha \goesto 0} \f{Dir}(\alpha, \v{F}_0)$ for $\v{F}_0$ arbitrary.
\\
$\f{DP}(\alpha, F_0)$ & The Dirichlet Process with concentration parameter $\alpha$ and mean CDF $F_0$.
\\
$\f{DP}(n, \hat{F}_n)$ & The Dirichlet Process whose mean CDF $\hat{F}_n$ is the empirical CDF of the data set $\v{y}$ of size $n$.
\\
$\f{DP}(0)$ & The improper Dirichlet Process corresponding to $\lim_{\alpha \goesto 0} \f{DP}(\alpha, \v{F}_0)$ for arbitrary $F_0$.
\\
\bottomrule
\end{tabular}
\end{center}
\caption{Notation. Bold symbols refer to vectors. Improper distributions are considered only as limits of conjugate families---we do not attempt to define $\f{DP}(0)$ directly as a non-normalizable measure. CDF: cumulative distribution function.}
\label{table-notation}
\end{table}

The $\f{DP}(n, \hat{F}_n)$ posterior and implied $\f{DP}(0)$ prior---see Table \ref{table-notation} for the notational conventions used in this work---have not been subject to the same level of formal study as Dirichlet Process Mixture priors and other priors over CDFs, in part due to the simplicity and discrete nature of $\f{DP}(n, \hat{F}_n)$.
On~the other hand, Dirichlet and Dirichlet Process priors with small concentration parameters have been used as low-information priors in a variety of settings (e.g., \cite{gelman14}), without much formal justification.
In this paper we offer a mathematical foundation showing that the use of $\f{DP}(0)$ is statistically sound.

\section{Results}

\subsection{Preliminaries}

To begin our discussion, we first introduce the notion of an invariant distribution, which describes what we mean by the term noninformative.

\begin{Definition}\label{def1}
[Invariant Distribution]
A density $\pi(\v{\theta})$ is \emph{invariant with respect to a transformation group} $G$ if for all $\tilde{\pi}(\v{\theta}) = \pi[g(\v{\theta})]$ with $g \in G$, and all measurable sets $A$,
\[ \label{invariant-0}
\int_A \pi(\v{\theta}) \dif\v{\theta} = \int_A \pi[g(\v{\theta})] \dif g(\v{\theta}) = \int_A \tilde{\pi} (\v{\theta}) \, \abs{\frac{\partial[g(\v{\theta})]}{\partial(\v{\theta})}} \dif \v{\theta} \, ,
\]
where $\abs{\frac{\partial[g(\v{\theta})]}{\partial(\v{\theta})}}$ is the Jacobian of the transformation.
\end{Definition}

Note that in Equation (\ref{invariant-0}), if we were to instead take $A$ in the middle and right integrals to be $g^{-1}(A)$, we would exactly get the classical integration by substitution formula, which under appropriate conditions is always true.~We are interested in the inverse problem: given a set of transformations in $G$, does there exist a unique $\pi$ satisfying (\ref{invariant-0})?

In a number of practically-relevant cases, $G$ is uniquely specified by the context of the problem being studied.
If this leads to a unique prior distribution $\pi$, and when additionally a unique likelihood also arises, for example via exchangeability, an optimal Bayesian analysis is possible, as defined in Section \ref{introduction}.
It is often the case that the prior distributions that result from this line of reasoning are limits of conjugate families, making them easy to work with---this occurs in our results below, in which the corresponding posterior distributions are Dirichlet.

The above definition is intuitive, but not sufficiently general to be applicable to spaces of functions.
There are multiple technical issues:
\begin{itemize}[leftmargin=2.3em,labelsep=4mm]
\item[(a)]  in many cases, $\pi$ cannot be taken to integrate to 1;
\item[(b)]  probability distributions on spaces of functions may not admit Riemann-integrable densities;
\item[(c)]  $G$ may be defined via equivalence classes of transformations, leading to singular Jacobians; and
\item[(d)]  infinite-dimensional measures that are non-normalizable are not well-behaved mathematically.
\end{itemize}

As a result, the above definition needs to be extended to a measure-theoretic setting.
We call a transformation group $G$ acting on a measure space \emph{nonsingular} if for $g \in G$ with $\tilde{\pi}(\v{\theta}) = \pi[g(\v{\theta})]$, we have $\pi \ll \tilde{\pi} \ll \pi$, where $\ll$ denotes absolute continuity of measures.

\begin{Definition}\label{def2}
[Invariant Measure]
Let $G$ be a nonsingular transformation group acting on a measure space.
We say that a measure $\pi$ is \emph{invariant with respect to $G$} if for any $g \in G$ with $\tilde{\pi}(\v{\theta}) = \pi[g(\v{\theta})]$ and for any measurable subset $A$ we have
\[ \label{invariant-1}
\int_\Omega \I_A \dif \pi = \int_\Omega \I_A \frac{\dif \tilde{\pi}}{\dif \pi} \dif\tilde{\pi} \, ,
\]
where $\Omega$ is the domain of $\pi$, $\I_A$ is the indicator function of the set $A$, and $\frac{\dif \tilde{\pi}}{\dif \pi}$ is the Radon–Nikodym  derivative of $\tilde{\pi}$ with respect to $\pi$.
\end{Definition}

It can be seen by taking $\pi$ to be absolutely continuous with respect to the Lebesgue measure that Equation (\ref{invariant-1}) is a direct extension of Equation (\ref{invariant-0}).

We would ultimately like to extend the above definition to the infinite-dimensional setting.
Doing so directly is challenging, because $\pi$ may be non-normalizable, in which case Kolmogorov's Consistency Theorem and other analytic tools for infinite-dimensional probability measures do not apply.
Here we sidestep this problem by instead extending the definition of invariance to allow us to define a sequence of \emph{approximately} invariant measures, which in our setting can be taken to be probability measures.
To do so, two additional definitions are needed.

\begin{Definition}[$\v{\eps}$-invariant Measure]
Let $G$ be a nonsingular transformation group acting on a measure space with invariant measure $\hat{\pi}$.
We say that a sequence of measures $\{\pi^{(\eps)} : \eps > 0\}$ is \emph{$\eps$-invariant with respect to $G$} if for any $g \in G$ with $\tilde{\pi}^{(\eps)}(\v{\theta}) = \pi^{(\eps)}[g(\v{\theta})]$ and each measurable subset $A$, the inequality
\[
\abs{\int_\Omega \I_A \dif \pi^{(\eps)} - \int_\Omega \I_A \frac{\dif \tilde{\pi}^{(\eps)}}{\dif \pi^{(\eps)}} \dif\tilde{\pi}^{(\eps)}} < \eps
\]
implies that
\[
\abs{\pi^{(\eps)}(A) - \hat{\pi}(A)} \leq \delta \mu(A) \, ,
\]
where $\mu(A)$ is a function, $\eps \goesto 0$ implies that $\delta \goesto 0$, and $\Omega$ is the domain of $\pi^{(\eps)}$ for all $\eps$.
\end{Definition}

\begin{Definition}[$\v{\eps}$-invariant Process] \label{epsilon-invariant-process}
Let $\{\Pi^{(\eps)} : \eps > 0\}$ be a sequence of stochastic processes, and let $G$ be a nonsingular transformation group.
Let $I$ be an arbitrary finite subset of the index set of the process, let $\pi^{(\eps)}_I$ be the finite-dimensional measure of $\Pi^{(\eps)}$ under $I$, and let $G_I$ be a finite-dimensional homomorphism of $G$ with invariant measure $\hat{\pi}_I$.
We say that the sequence of processes $\Pi^{(\eps)}$ is \emph{$\eps$-invariant} if, for each $I$, each $g_I \in G_I$ with $\tilde{\pi}^{(\eps)}_I(\v{\theta}) = \pi^{(\eps)}_I[g_I(\v{\theta})]$ and each measurable subset $A$, the inequality
\[
\abs{\int_{\Omega_I} \I_A \dif \pi^{(\eps)}_I - \int_{\Omega_I} \I_A \frac{\dif \tilde{\pi}^{(\eps)}_I}{\dif \pi^{(\eps)}_I} \dif\tilde{\pi}^{(\eps)}_I} < \eps
\]
implies that
\[
\abs{\pi^{(\eps)}_I(A) - \hat{\pi}_I(A)} \leq \delta \mu_I(A) \, ,
\]
where $\mu_I(A)$ is a function, $\eps \goesto 0$ implies that $\delta \goesto 0$, $\Omega_I$ is the domain of $\pi^{(\eps)}_I$ for all $\eps$, and $( \eps, \delta )$ can be taken to be identical for all $I$.
\end{Definition}

Definition \ref{epsilon-invariant-process} has been explicitly chosen to formalize the notion of noninformativeness on a space of functions without constructing a non-normalizable infinite-dimensional measure.

To \textls[-15]{complete our assumptions, we need to specify $G$.
Our definitions constitute a direct generalization of the transformation group used by Jaynes to derive the Haldane prior for $p=2$---see Section} \ref{introduction}.

\begin{Definition}[Probability Function Transformation Group] \label{group-functions}
Let
\[
G_\infty = \cbr[1]{g \! : S_\infty \goesto S_\infty}
\]
be a\textls[-15]{ nonsingular group of measurable functions under composition acting on the infinite-dimensional simplex} $S_\infty$.
\end{Definition}

\begin{Definition}[Probability Vector Transformation Group] \label{group-vectors}
For non-negative integer $p$ and any vector $( c_1, .. , c_p )$ of non-negative constants, let
\[ \label{basic-transformation-1}
G_p = \cbr{g \! : (\theta_1, .., \theta_p) \goesto \del{\frac{c_1 \theta_1}{\sum_{i=1}^p c_i \theta_i}, .., \frac{c_p \theta_p}{\sum_{i=1}^p c_i \theta_i}}}
\]
be a nonsingular group under composition acting on the $p$-dimensional simplex $S_p$, where each element $g \in G$ represents an equivalence class of the transformations (\ref{basic-transformation-1}).
\end{Definition}

Note that $G_p$ is a $p$-dimensional homomorphism of $G_\infty$---we use this property in our proofs below.
It can also readily be seen that for any $g$, the constants $c_i$ are only determined up to proportionality.

\begin{Proposition}[Radon–Nikodym Derivative]
For each $g \in G_p$ and $\tilde{\pi}(\v{\theta}) = \pi[g(\v{\theta})]$, the Radon–Nikodym derivative of $\tilde{\pi}$ with respect to $\pi$ is
\[
\frac{\dif \tilde{\pi}}{\dif \pi} (\v{\theta}) = \frac{\prod_{i=1}^p c_i}{\del{\sum_{i=1}^p c_i \theta_i}^p}
.
\]
\end{Proposition}

\begin{proof}
Let $\lambda$ be the Lebesgue measure on the $p$-dimensional probability simplex, and define \mbox{$\tilde{\lambda}(\v{\theta}) = \lambda[g(\v{\theta})]$}.
Note first that $\lambda \ll \tilde{\lambda} \ll \pi \ll \tilde{\pi} \ll \lambda$.
Note also that
\[
\frac{\dif \pi}{\dif \lambda} = \frac{\dif \tilde{\pi}}{\dif \tilde{\lambda}} \, ,
\]
because the same transformation $g$ is used in defining $\tilde{\pi}$ and $\tilde{\lambda}$.
Then, note that
\[
\frac{\dif \tilde{\pi}}{\dif \pi} = \frac{\dif \tilde{\pi}}{\dif \pi} \frac{\dif \tilde{\lambda}}{\dif \lambda} \frac{\dif \lambda}{\dif \tilde{\lambda}} = \frac{\dif \lambda}{\dif \pi} \frac{\dif \tilde{\lambda}}{\dif \lambda} \frac{\dif \tilde{\pi}}{\dif \tilde{\lambda}} = \frac{\dif \tilde{\lambda}}{\dif \lambda} \, ,
\]
and hence it suffices to consider the transformation $g$ applied to the Lebesgue measure.
Consider an arbitrary hypercube $B$.
We have
\[
\lambda(B) = \lambda_1 (B_1) \, .. \, \lambda_p (B_p) \, ,
\]
where $\lambda_i$ are 1-dimensional Lebesgue measures, for which we have that
\[
\lambda_i (B_i) = b_i - a_i \, ,
\]
where $[a_i, b_i]$ is the one-dimensional projection of the hypercube $B$ in dimension $i$.
Consider now the transformation $g$.
We may decompose $g$ into $d$ and $n$ where
\<
d &: (\theta_1,..,\theta_p) \goesto (c_1 \theta_1,..,c_p \theta_p)
&
n &: (\theta_1,..,\theta_p) \goesto \del{\frac{\theta_1}{\sum_{i=1}^p \theta_i},..,\frac{\theta_p}{\sum_{i=1}^p \theta_p}}
\,.
\>

Now consider the effect of $d$ and $n$ on $\lambda_i$.
We have
\<
\lambda_i [d(B_i)] &= c_i(b_i - a_i)
&
&\text{and}
&
\lambda_i [n(B_i)] &= \frac{b_i - a_i}{\sum_{i=1}^p (b_i - a_i)} \,,
\>
hence
\[
\lambda_i [g(B_i)] = \frac{c_i(b_i - a_i)}{\sum_{j=1}^p c_j(b_j - a_j)}
\,.
\]

Therefore
\[
\lambda[g(B)] = \prod_{i=1}^p \frac{c_i(b_i - a_i)}{\sum_{j=1}^p c_j(b_j - a_j)}
\]
and we can compute the ratio
\[
\frac{\tilde{\lambda}(B)}{\lambda(B)} = \frac{\lambda[g(B)]}{\lambda(B)} = \prod_{i=1}^p \frac{c_i(b_i - a_i)}{\sum_{j=1}^p c_j(b_j - a_j)} \sbr{\prod_{i=1}^p (b_i - a_i)}^{-1} = \frac{\prod_{i=1}^p c_i}{\sbr{\sum_{i=1}^p c_i(b_i - a_i)}^p}
\,.
\]

This holds for all $B$, hence the Radon–Nikodym derivative is just
\[
\frac{\dif \tilde{\lambda}}{\dif \lambda} (\v{\theta}) = \frac{\prod_{i=1}^p c_i}{\del{\sum_{i=1}^p c_i \theta_i}^p} \, ,
\]
which is the desired result.
\end{proof}

Since we \textls[-15]{are working with non-normalizable measures as improper priors, we cannot rigorously talk about their probability densities.
In many cases, such improper priors can be shown to be limits of families of conjugate priors for which the limiting posterior distribution is well-defined, making them usable in practice.
To make our discussion of improper priors rigorous, we need the following} definition.

\begin{Definition}[Generalized Density]
Let $\pi$ be a measure on $\R^p$ (for $p$ a positive integer) such that \mbox{$\pi \ll \lambda \ll \pi$}, where $\lambda$ is Lebesgue measure on $\R^p$.
Suppose that the Radon–Nikodym derivative of $\pi$ with respect to $\lambda$ is Riemann-integrable, and define a family of functions equal to the Radon–Nikodym derivative up to a proportionality constant.
We call any function in this family a \emph{generalized density} of $\pi$.
\end{Definition}

\subsection{Main Results}

\begin{Remark}[Notation] \label{notation}
In the following results, we will assume that $(\theta_1,..,\theta_p)$ is a probability vector of dimension $p \geq 2$.
$G_\infty$ and $G_p$ will be the transformation groups identified in Definitions \ref{group-functions} and \ref{group-vectors}, respectively.
As noted previously in Table \ref{table-notation}, $\f{Dir}(\alpha, \v{F}_0)$ will denote the Dirichlet distribution under the alternative parametrization based on concentration parameter $\alpha$ and mean probability vector $\v{F}_0$.
This is equivalent to the usual parameterization in terms of concentration vector $\v{\alpha}$ by the identity $\v{\alpha} = \alpha \v{F}_0$---we refer to this as the $\f{Dir}(\v{\alpha})$ distribution.
Similarly, $\f{DP}(\alpha, F_0)$ will refer to the Dirichlet Process with concentration parameter $\alpha$ and mean function $F_0$.
We will refer to the improper priors defined via the conjugate limits as $\alpha \goesto 0$ of $\f{Dir}(\alpha, F_0)$ and $\f{DP}(\alpha, F_0)$ for arbitrary $F_0$ as $\f{Dir}(0)$ and $\f{DP}(0)$, respectively.
\end{Remark}

We are now ready to introduce our first result.
The argument below is a direct generalization of the line of reasoning in \textcite{jaynes68}: the Haldane prior obtained is a special case of our result for $p=2$.

\begin{Theorem} \label{epsilon-invariance}
Among \textls[-25]{the class of measures that admit generalized densities, the measure $\pi$ with generalized} density
\[
\pi\del{\theta_1,..,\theta_p} \propto \frac{1}{\prod_{i=1}^p \theta_i} \,,
\]
which we call $\f{Dir}(0)$, is the unique invariant measure under $G_p$.
\end{Theorem}

\begin{proof}
An invariant measure $\pi$ under $G_p$ needs to satisfy the equation
\[ \label{theorem-10-1}
\int_{S_p} \I_A \dif \pi = \int_{S_p} \I_A \frac{\dif \tilde{\pi}}{\dif \pi} \dif\tilde{\pi} \, ,
\]
where $S_p$ is the $p$-dimensional simplex and $\tilde{ \pi } ( \v{ \theta } ) = \pi [ g ( \v{ \theta } ) ]$ for some $g \in G_p$.
Since $\pi$ is assumed to admit a generalized density, we can rewrite (\ref{theorem-10-1}) as a Riemann integral.
In addition, we substitute in the transformation and Radon–Nikodym derivative, and get
\[
\int_A \pi\del{\theta_1,..,\theta_p} \dif\theta_1 .. \dif\theta_p = \int_A \pi\del{\frac{c_1 \theta_1}{\sum_{i=1}^p c_i \theta_i}, .., \frac{c_p \theta_p}{\sum_{i=1}^p c_i \theta_i}} \frac{\prod_{i=1}^p c_i}{\del{\sum_{i=1}^p c_i\theta_i}^p} \dif\theta_1 .. \dif\theta_p
\, .
\]

This formula needs to hold for all measurable sets $A$, and hence the functions inside the integrals need to be equal pointwise.
This yields the functional equation

\[ \label{functional-equation-1}
\pi\del{\theta_1,..,\theta_p} = \pi\del{\frac{c_1 \theta_1}{\sum_{i=1}^p c_i \theta_i}, .., \frac{c_p \theta_p}{\sum_{i=1}^p c_i \theta_i}} \frac{\prod_{i=1}^p c_i}{\del{\sum_{i=1}^p c_i\theta_i}^p} \, ,
\]
which will be the main subject of further study.
This is a multivariate functional equation that at first may appear fearsome, but is in fact solvable via elementary methods.
To solve it, recognizing that (\ref{functional-equation-1}) must hold for all probability vectors $(\theta_1,..,\theta_p)$ and all vectors $( c_1, .., c_p )$ of positive constants $c_i$, we set
\<
\del{\theta_1,..,\theta_p} &= \del{p^{-1},..,p^{-1}}
&
&\text{and}
&
\sum_{i=1}^p c_i = 1 \, ,
\>
which yields
\[ \label{functional-equation-2}
\pi\del{p^{-1},..,p^{-1}} = \pi\del{\frac{c_1 p^{-1}}{ p^{-1} \sum_{i=1}^p c_i}, .., \frac{c_p p^{-1}}{p^{-1}\sum_{i=1}^p c_i}} \frac{\prod_{i=1}^p c_i}{\del{p^{-1}\sum_{i=1}^p c_i}^p} \, .
\]

Then, by swapping $c_i$ for $\theta_i$, (\ref{functional-equation-2}) rearranges into
\[ \label{functional-equation-3}
\pi\del{\theta_1,..,\theta_p} = \frac{\pi\del{p^{-1},..,p^{-1}}p^{-p}}{\prod_{i=1}^p \theta_i} \propto \frac{1}{\prod_{i=1}^p \theta_i} \, ,
\]
since the numerator is not a function of any $\theta_i$, and it can easily be checked that all such generalized densities are valid solutions to the original equation.
Thus (\ref{functional-equation-3}) is the functional equation's unique solution and therefore the unique invariant measure under $G_p$.
\end{proof}

The same technique used to solve the functional equation in Theorem \ref{epsilon-invariance} can be used to prove a much stronger result: if the functional equation is true approximately, its solutions will approximate those of the exact equation.
In the next result we make use of the definition of \emph{stability} of a functional equation due to Hyers, Ulam and Rassias---see \textcite{jung11} for details.

\begin{Corollary}[Hyers–Ulam–Rassias Stability] \label{stability}
Suppose we have
\[
\abs{\pi\del{\theta_1,..,\theta_p} - \pi\del{\frac{c_1 \theta_1}{\sum_{i=1}^p c_i \theta_i}, .., \frac{c_p \theta_p}{\sum_{i=1}^p c_i \theta_i}} \frac{\prod_{i=1}^p c_i}{\del{\sum_{i=1}^p c_i\theta_i}^p}} < \delta
.
\]

Then
\<
\abs{\pi\del{\theta_1,..,\theta_p} - \hat{\pi}\del{\theta_1,..,\theta_p} } &< \delta \frac{e^{e^{-1}}} {\prod_{i=1}^p \theta_i} \, ,
&
&\text{where}
&
\hat{\pi}\del{\theta_1,..,\theta_p} &\propto \frac{1}{\prod_{i=1}^p \theta_i}
.
\>
\end{Corollary}

\begin{proof}
By repeating the technique from the previous proof, we have
\[
\abs{\pi\del{p^{-1},..,p^{-1}} - \pi\del{c_1,..,c_p} \frac{\prod_{i=1}^p c_i}{p^{-p}}} < \delta \, ,
\]
which can be rewritten
\[
\abs{\pi\del{\theta_1,..,\theta_p} - \frac{\pi\del{p^{-1},..,p^{-1}} p^{-p}}{\prod_{i=1}^p \theta_i}} < \delta \frac{p^{-p}} {\prod_{i=1}^p \theta_i} < \delta \frac{e^{e^{-1}}} {\prod_{i=1}^p \theta_i} \, ,
\]
where the last inequality is strict because $p$ is a positive integer.
Letting
\[
\frac{\pi\del{p^{-1},..,p^{-1}}p^{-p}}{\prod_{i=1}^p \theta_i} \propto \frac{1}{\prod_{i=1}^p \theta_i} \propto \hat{\pi}\del{\theta_1,..,\theta_p} \, ,
\]
we get
\[
\abs{\pi\del{\theta_1,..,\theta_p} - \hat{\pi}\del{\theta_1,..,\theta_p} } < \delta \frac{e^{e^{-1}}} {\prod_{i=1}^p \theta_i} \, ,
\]
which is the stability result desired.
\end{proof}

This suffices to prove our result for the Dirichlet distribution.

\begin{Theorem} \label{dirichlet-epsilon-invariance}
$\f{Dir}(\eps, \v{F}_0)$ is an $\eps$-invariant measure under $G_p$ for all $\v{F}_0$.
\end{Theorem}

\begin{proof}
By repeating the steps of Theorem \ref{epsilon-invariance} and combining them with Corollary \ref{stability}, we obtain that $\f{Dir}(\eps, \v{F}_0)$ is $\eps$-invariant under $G_p$ if and only if it satisfies
\<
\abs{\pi^{(\eps)}\del{\theta_1,..,\theta_p} - \hat{\pi}\del{\theta_1,..,\theta_p} } &< \delta \frac{e^{e^{-1}}} {\prod_{i=1}^p \theta_i}
&
&\text{for some}
&
\hat{\pi}\del{\theta_1,..,\theta_p} &\propto \frac{1}{\prod_{i=1}^p \theta_i}
.
\>

Substituting in $\f{Dir}(\eps, \v{F}_0)$, and choosing the constant $C_\eps$ of the generalized density $\hat{\pi}$ to be the same as for the Dirichlet, we get
\[
\abs{C_\eps \prod_{i=1}^p \theta_i^{\, \eps F_{0i} - 1} - \frac{C_\eps}{\prod_{i=1}^p \theta_i} } < \delta \frac{e^{e^{-1}}} {\prod_{i=1}^p \theta_i} \, ,
\]
where $F_{0i}$ are the components of the probability vector $\v{F}_0$, and this expression simplifies to
\[
C_\eps \, e^e \, \abs{\prod_{i=1}^p \theta_i^{\, \eps F_{0i}} - 1 } < \delta
.
\]

Since $0 \leq \theta_i \leq 1$ for all $i$, the product is upper bounded by $1$ and lower bounded by $0$.
Thus the inequality holds near zero if
\[
C_\eps < \delta
\]
for all $( \theta_1,..,\theta_p )$, and since $C_\eps \goesto 0$ we get that, as $\eps \goesto 0$, we can choose $\delta$ such that $\delta \goesto 0$.
Thus,~$\f{Dir}(\eps, \v{F}_0)$ is $\eps$-invariant for all $\v{F}_0$.
\end{proof}

We now extend Theorem \ref{dirichlet-epsilon-invariance} to get an analogous result for the Dirichlet Process.

\begin{Theorem} \label{dp-epsilon-invariance}
$\f{DP}(\eps, F_0)$ is an $\eps$-invariant process under $G_\infty$ for all $F_0$.
\end{Theorem}

\begin{proof}
Consider an arbitrary finite-dimensional index $I$ with corresponding homomorphism $G_I$ and finite-dimensional measure $\pi^{(\eps)}_I$.
It follows from Theorem \ref{dirichlet-epsilon-invariance} that $\pi^{(\eps)}_I$ is $\eps$-invariant with
\[
C_\eps < \delta
.
\]

This inequality depends only on $C_\eps$, so it suffices to show that this constant can be bounded by another constant that is not a function of $p$ and approaches 0.
$C_\eps$ is an instance of the inverse multivariate beta function, which is a ratio of gamma functions.
It is well known that
\[
\lim_{x \goesto 0} \sbr{ \frac{1}{x} - \Gamma(x) } = \gamma ,
\]
where $\gamma$ is the Euler-Mascheroni constant.
Therefore, we have
\[
C_\eps = \frac{\Gamma(\eps)}{\prod_{i=1}^p \Gamma(\eps F_{0i})} = \frac{O(1/\eps)}{\prod_{i=1}^p O(1/\eps) } \leq \frac{O(1/\eps)}{\prod_{i=1}^2 O(1/\eps)} = O(\eps) \goesto 0
\]
as $\eps \goesto 0$.
Thus, for each $\eps$, we can choose a $\delta$ to satisfy the required expressions under all finite-dimensional index sets, and $\f{DP}(\eps, F_0)$ is therefore an $\eps$-invariant process.
\end{proof}

\scalebox{0.95}[1]{We conclude our theoretical investigation with a conjecture: the $\eps$-invariance of all finite-dimensional} distributions with a uniform $\delta$ should suffice for invariance with respect to the original group acting on the infinite-dimensional space.

\begin{conjecture} \label{conjecture}
A stochastic process is an $\eps$-invariant process if and only if the measure of its sample paths is an $\eps$-invariant measure.
\end{conjecture}

One approach to attempting a proof of this conjecture would involve appropriately extending Kolmogorov's Consistency Theorem to $\sigma$-finite infinite-dimensional measures.
This can be done, but the notions involved are quite technical---see \textcite{yamasaki85} for more details.

\section{Discussion} \label{discussion}

To see how our results may be applied, consider again the randomized controlled trial of Section~\ref{introduction}, and suppose now that the outcome $y_i$ for participant $i$ in the experimental group is categorical with $p$ levels.
Under exchangeability, a minor extension of de Finetti's Theorem for dichotomous outcomes then yields that the likelihood can be expressed as
\[
y_i \given \v{\theta} \iid[MN](1, \v{\theta}) \, ,
\]
in which MN$(k, \v{\theta} )$ is the multinomial distribution with parameters $k$ and $\v{ \theta }$.
Theorem \ref{epsilon-invariance} implies that, modulo inherent abuse of notation under improper priors,
\[
(\theta_1,..,\theta_p) \avtdist[Dir](0)
\]
is the unique prior that obeys the fundamental invariance possessed by the problem---namely, invariance with respect to all transformations of probability vectors that preserve normalization.
Thus we have extended Jaynes's result for binomial outcomes to the multinomial setting, yielding another instance of optimal Bayesian analysis.

Generalizing to the setting where $\v{y}$ is an exchangeable sequence of real-valued outcomes, de~Finetti's most general representation theorem implies that
\[
y_i \given F \iid F
\]
is the unique likelihood.
If little is known about $F$, and it is therefore approximately invariant under all measurable functions---i.e., under $G_\infty$, see Definition \ref{group-functions}---the prior given by Theorem \ref{dp-epsilon-invariance} is
\[ \label{dp-eps-F0}
F \avtdist[DP](\eps, F_0)
\,.
\]

By the usual conjugate updating in the Dirichlet Process setting, the posterior on $F$ given $\v{ y }$ with the prior in (\ref{dp-eps-F0}) is
\[ \label{de-finetti-representation-3}
F \given \v{y} \avtdist[DP] \del{\eps + n, \frac{\eps}{\eps + n} F_0 + \frac{n}{\eps+n} \hat{F}_n} \, ,
\]
in which $\hat{ F }_n$ is the empirical CDF based on $\v{ y }$. Since $\eps$ may be taken as close to zero as one wishes, it is natural to regard

\[
F \given \v{y} \avtdist[DP]\del{n, \hat{F}_n}
\]
as an instance of approximately optimal Bayesian analysis for all $\eps$.
Conjecture \ref{conjecture} would strengthen this assertion---provided $\f{DP}(0)$ can be rigorously constructed as an infinite-dimensional $\sigma$-finite measure, which is beyond the scope of this work.

Though the simplicity of this analysis may at first make it seem limited, its appeal comes from its extremely general ability to characterize uncertainty.
See, e.g., \textcite{terenin15} for an example of a $\f{DP}(n, \hat{F}_n)$ analysis in two randomized controlled trials in e-commerce, one with sample sizes in the tens of millions.
Furthermore, sampling from $\f{DP}(n, \hat{F}_n)$ on a discrete domain has recently been shown in a completely different setting---see Appendix B of Terenin et al. \cite{terenin17b}---to be asymptotically equivalent to the widely-used frequentist bootstrap of \textcite{efron79}. This also applies to the Bayesian bootstrap of \textcite{rubin81}, since it is asymptotically equivalent to the frequentist version.
Our analysis provides a Bayesian nonparametric justification for this class of methods.

Bayesian analysis cannot proceed without the specification of a stochastic model---prior and sampling distribution---relating known quantities to unknown quantities: data to parameters.
One of the great challenges of applied statistics is that the model is not necessarily uniquely determined by the context of the problem under study, giving rise to model uncertainty, which if not assessed and correctly propagated can cause badly calibrated and unreliable inference, prediction and decision---see, e.g., \textcite{draper95}.
Perhaps the simplest way to avoid model uncertainty is to recognize settings in which it does not exist---situations where broad and simple mathematical assumptions, rendered true by problem context, lead to unique posterior distributions. Our term for this is \emph{optimal Bayesian analysis}.
It seems worthwhile (a) to catalog situations in which optimal analysis is possible and (b) to work to extend the list of such situations---Theorems \ref{epsilon-invariance} and \ref{dp-epsilon-invariance} are two contributions to this effort.

\acknowledgments{We are grateful to Daniele Venturi, Yuanran Zhu, and Catherine Brennan for their thoughts on differential equations, which we originally used in a much longer and more complicated proof of the solution of our functional equation.
We are grateful to Dan Simpson for his thoughts on infinite-dimensional measures.
We are additionally grateful to Juhee Lee for her comments on prior specification, and to Thanasis Kottas for his thoughts on Dirichlet Processes.
Membership on this list does not imply agreement with the ideas expressed here, nor are any of these people responsible for any errors that may be present.}
\authorcontributions{Alexander Terenin and David Draper contributed to the conceptual and theoretical development of the methods in this work and co-wrote the manuscript.}
\conflictsofinterest{The authors declare no conflict of interest.}


\begin{thebibliography}{999}
\providecommand{\natexlab}[1]{#1}

\bibitem[Fisher(1925)]{fisher25}
Fisher, R.A.
\newblock {\em Statistical Methods for Research Workers}; Oliver and Boyd: Edinburgh, UK,
  1925.

\bibitem[{Medical Research Council}(1948)]{mrc48}
{Medical Research Council}.
\newblock Streptomycin treatment of pulmonary tuberculosis.
\newblock {\em Br. Med. J.} {\bf 1948}, {\em 2},~769--782.

\bibitem[Kohavi and Longbotham(2015)]{longbotham15}
Kohavi, R.; Longbotham, R.
\newblock Online controlled experiments and AB tests. In {\em Encyclopedia of
  Machine Learning and Data Mining}; Springer:  Berlin, Germany, 2015.

\bibitem[de~Finetti(1937)]{definetti37}
De~Finetti, B.
\newblock La pr{\'e}vision: Ses lois logiques, ses sources subjectives.
\newblock {\em Annales de l'institut Henri Poincar{\'e}} {\bf 1937}, {\em
  7},~1--68. (In French).

\bibitem[Jeffreys(1946)]{jeffreys46}
Jeffreys, H.
\newblock An invariant form for the prior probability in estimation problems.
\newblock  \emph{Proc. R. Soc. Lond. A Math. Phys. Eng. Sci.}  \textbf{1946}, \emph{186}, 453--461.

\bibitem[Jaynes(2003)]{jaynes03}
Jaynes, E.T.
\newblock {\em Probability Theory: The Logic of Science}; Cambridge University Press: Cambridge, UK,  2003.

\bibitem[Eaton(1989)]{eaton89}
Eaton, M.L.
\newblock {\em Group Invariance Applications in Statistics}; Regional
  Conference Series in Probability and Statistics; Institute of Mathematical Statistics: Shaker Heights, OH, USA, 1989.

\bibitem[Ferguson(1973)]{ferguson73}
Ferguson, T.S.
\newblock A Bayesian analysis of some nonparametric problems.
\newblock {\em Ann. Stat.} {\bf 1973}, \emph{1}, 209--230.

\bibitem[Laplace(1774)]{laplace74}
Laplace, P.S.
\newblock M\'emoire sur la probabilit\'e des causes par les \'ev\'enements.
\newblock {\em M\'emoires de l'Acad\'emie Royale des Sciences de Paris} {\bf
  1774}, {\em 6},~621. (In French).

\bibitem[Hald(2007)]{hald07}
Hald, A.
\newblock {\em A History of Parametric Statistical Inference from Bernoulli to
  Fisher, 1713--1935}; Springer: Berlin, Germany, 2007.

\bibitem[Bernardo(1979)]{bernardo79}
Bernardo, J.M.
\newblock Reference posterior distributions for Bayesian inference.
\newblock {\em J. R. Stat. Soc. Ser. B
  (Methodol.)} {\bf 1979}, \emph{41}, 113--147.

\bibitem[Jaynes(1968)]{jaynes68}
Jaynes, E.T.
\newblock Prior probabilities.
\newblock {\em IEEE Trans. Syst. Sci. Cybern.} {\bf
  1968}, {\em 4},~227--241.

\bibitem[Bush \em{et~al.}(2010)Bush, Lee, and MacEachern]{bush10}
Bush, C.A.; Lee, J.; MacEachern, S.N.
\newblock Minimally informative prior distributions for non-parametric Bayesian
  analysis.
\newblock {\em J. R. Stat. Soc. Ser. B (Stat.
  Methodol.)} {\bf 2010}, {\em 72},~253--268.

\bibitem[Lee \em{et~al.}(2014)Lee, MacEachern, Lu, and Mills]{lee14}
Lee, J.; MacEachern, S.N.; Lu, Y.; Mills, G.B.
\newblock Local-mass preserving prior distributions for nonparametric Bayesian
  models.
\newblock {\em Bayesian Anal.} {\bf 2014}, {\em 9},~307--330.

\bibitem[Gelman \em{et~al.}(2014)Gelman, Carlin, Stern, Dunson, Vehtari, and
  Rubin]{gelman14}
Gelman, A.; Carlin, J.B.; Stern, H.S.; Dunson, D.B.; Vehtari, A.; Rubin, D.B.
\newblock {\em Bayesian Data Analysis}, 3rd ed.; CRC~Press: Boca Raton, FL, USA,  2014.

\bibitem[Jung(2011)]{jung11}
Jung, S.M.
\newblock {\em Hyers-Ulam-Rassias Stability of Functional Equations in
  Nonlinear Analysis}; Springer: Berlin, Germany,  2011.

\bibitem[Yamasaki(1985)]{yamasaki85}
Yamasaki, Y.
\newblock {\em Measures on Infinite-Dimensional Spaces}; World Scientific: Singapore,
  1985.

\bibitem[Terenin and Draper(2015)]{terenin15}
\scalebox{0.92}[1]{Terenin, A.; Draper, D.
 Cox's Theorem and the Jaynesian Interpretation of Probability.
{\em arXiv} {\bf 2015}, arXiv:1507.06597.}

\bibitem[Terenin \em{et~al.}(2017)Terenin, Magnusson, Jonsson, and
  Draper]{terenin17b}
Terenin, A.; Magnusson, M.; Jonsson, L.; Draper, D.
\newblock P{\'o}lya Urn Latent Dirichlet Allocation: A doubly sparse massively
  parallel sampler.
\newblock {\em arXiv} {\bf 2017}, arXiv:1704.03581.

\bibitem[Efron(1979)]{efron79}
Efron, B.
\newblock Bootstrap methods: Another look at the jackknife.
\newblock {\em Ann. Stat.} {\bf 1979}, {\em 7},~1--26.

\bibitem[Rubin(1981)]{rubin81}
Rubin, D.B.
\newblock The Bayesian Bootstrap.
\newblock {\em Ann. Stat.} {\bf 1981}, {\em 9},~130--134.

\bibitem[Draper(1995)]{draper95}
Draper, D.
\newblock Assessment and propagation of model uncertainty.
\newblock {\em J. R. Stat. Soc. Ser. B (Stat.
  Methodol.)} {\bf 1995}, {\em 57}, 45--97.

\end{thebibliography}
\end{document}